\documentclass[12pt]{article}
\usepackage{amscd,amssymb,epic,amsmath}
\usepackage[all]{xy}

\title{Regular functions on the Shilov boundary}\author{Olga Bershtein}

\begin{document}
\newtheorem{proposition}{Proposition}
\newtheorem{lemma}{Lemma}
\newtheorem{corollary}{Corollary}

\maketitle
\centerline{email:bershtein@ilt.kharkov.ua}
\begin{abstract}
In this paper a $*$-algebra of regular functions on the Shilov boundary
$S(\mathbb D)$ of bounded symmetric domain $\mathbb D$ is constructed. The
algebras of regular functions on $S(\mathbb D)$ are described in terms of
generators and relations for two particular series of bounded symmetric
domains. Also, the degenerate principal series of quantum Harish-Chandra
modules related to $S(\mathbb D)=U_n$ is investigated.
\end{abstract}

\section{Introduction} Quantum Harish-Chandra modules (see the next Sec.)
form a broad class of infinite dimensional representations of quantum
universal enveloping algebras. Construction and classification of simple
quantum Harish-Chandra modules are important open problems (cf. \cite{Shm}).

In this paper we investigate a quantum analog of degenerate principal
series realized in the spaces of regular functions on the Shilov boundary
of bounded symmetric domain of tube type \cite{Wolf}.

Let $(a_{ij})_{i,j=1,\ldots,l}$ be a Cartan matrix of positive type,
$\mathfrak g$ the corresponding simple complex Lie algebra. So the Lie
algebra can be defined by the generators $e_i,f_i,h_i,i=1,...,l$, and the
well-known relations (see \cite{Jantzen}). Let $\mathfrak{h}$ be the
linear span of $h_i,i=1,...,l$. Fix simple roots $\{\alpha_i \in
\mathfrak{h}^*|i=1,...,l\}$ via $\alpha_i(h_j)=a_{ji}$. Also, let
$\{\varpi_i| i=1,...,l\}$ be the fundamental weights and
$P=\bigoplus_{i=1}^l \mathbb{Z}\varpi_i$ the weight lattice.

Fix $l_0 \in \{1,...,l\}$ and the Lie subalgebra $\mathfrak{k} \subset
\mathfrak{g}$ generated by
$$
e_i, f_i,\;\; i\neq l_0; \qquad h_i, \;\;i=1,...,l.
$$
We consider Lie algebras $\mathfrak{g}$ equipped with $\mathbb Z$-grading
as follows:
\begin{equation}\label{par_type}
\mathfrak{g}=\mathfrak{g}_{-1}\oplus\mathfrak{g}_0\oplus\mathfrak{g}_{+1},
\quad  \mathfrak{g}_j=\{\xi\in\mathfrak{g}|\:[h_0,\xi]=2j\xi\},
\end{equation}
where $h_{0} \in \mathfrak{h}$ and
$$
\alpha_i(h_0)=0, \; i \neq l_0; \qquad \alpha_{l_0}(h_0)=2.
$$
Let $\delta$ be the maximal root, and $\delta=\sum_{i=1}^l c_i\alpha_i$.
(\ref{par_type}) holds if and only if $c_{l_0}=1$. In this case
$\mathfrak{g}_0=\mathfrak{k}$ and we call the pair
$(\mathfrak{g},\mathfrak{k})$ a Hermitian symmetric pair.

Harish-Chandra established the one-to-one correspondence between Hermitian
symmetric pairs and irreducible bounded symmetric domains considered up to
biholomorphic isomorphisms \cite{Wolf}.

Let $W$ be the Weyl group of the root system $R$, $w_0 \in W$ the element of
maximal length. The irreducible bounded symmetric domain $\mathbb D$ related
to the pair $(\mathfrak{g},\mathfrak{k})$ is a domain of tube type if and
only if $\varpi_{l_0}=-w_0\varpi_{l_0}$.

\medskip

Fix a Hermitian symmetric pair $(\mathfrak{g},\mathfrak{k})$. Let $G$ be a
connected simply connected complex Lie group with ${\rm Lie}\, (G)
=\mathfrak{g}$ and $K \subset G$ a connected complex Lie subgroup with
${\rm Lie}\,(K) =\mathfrak{k}$. Consider the Lie subalgebra $\mathfrak{p}
\subset \mathfrak{g}$ generated by $e_i, h_i, i =1,...,l$, and $f_j, j
\neq l_0$. There is the corresponding connected complex Lie subgroup $P
\subset G$. The homogeneous space $G/P$ is a simply connected projective
variety.

There exists a distinguished noncompact real form $G_0$ of $G$ with a unique
closed $G_0$-orbit in $G/P$ \cite{Wolf}. The Shilov boundary $S(\mathbb D)$
corresponds to it under the Borel embedding $i:\mathbb D \hookrightarrow
G/P$. These reasons allow us to obtain quantum analogs of regular functions
on the Shilov boundary and series of Harish-Chandra modules related to it.
This approach to the quantum Shilov boundaries belongs to L. Vaksman
(private communication).

\section{Quantum analog of algebra of regular functions on the Shilov
boundary}

In this section we introduce quantum analogs of regular functions on the
Shilov boundaries $S(\mathbb D)$ of bounded symmetric domains $\mathbb D$
of tube type.

First of all, recall some notions from the quantum group theory
\cite{Jantzen}. In the sequel the ground field is $\mathbb C$, $q \in
(0,1)$, and all algebras are associative and unital.

Denote by $d_i>0, i=1,\ldots,l$, such coprime numbers that the matrix
$(d_ia_{ij})_{i,j=1,\ldots,l}$ is symmetric. Recall that the quantum
universal enveloping algebra $U_q \mathfrak{g}$ is a Hopf algebra defined
by the generators $K_i$, $K_i^{-1}$, $E_i$, $F_i$, $i=1,2,\ldots,l$ and
the relations:
$$K_iK_j=K_jK_i,\quad K_iK_i^{-1}=K_i^{-1}K_i=1,$$
\begin{equation*}
K_iE_j=q_i^{a_{ij}}E_jK_i,\quad K_iF_j=q_i^{-a_{ij}}F_jK_i,
\end{equation*}
$$
E_iF_j-F_jE_i=\delta_{ij}\,\frac{K_i-K_i^{-1}}{q_i-q_i^{-1}},
$$
\begin{equation*}
  \sum\limits_{m=0}^{1-a_{ij}}(-1)^m\left[1-a_{ij} \atop m\right]_{q_i}
  E_i^{1-a_{ij}-m}E_jE_i^m=0,
\end{equation*}
\begin{equation*}
  \sum\limits_{m=0}^{1-a_{ij}}(-1)^m\left[1-a_{ij} \atop m\right]_{q_i}
  F_i^{1-a_{ij}-m}F_jF_i^m=0,
\end{equation*}
where $q_i=q^{d_i}$, $1\leq i\leq l$ and
\begin{equation*}
  \left[m \atop n\right]_q=\frac{[m]_q!}{[n]_q![m-n]_q!},\quad
  [n]_q!=[n]_q\ldots[2]_q[1]_q,\quad
  [n]_q=\frac{q^n-q^{-n}}{q-q^{-1}}.
\end{equation*}

The comultiplication $\Delta$, the counit $\varepsilon$, and the antipod
$S$ are defined on generators by the following formulas:
\begin{equation*}
  \Delta(E_i)=E_i\otimes 1+K_i\otimes E_i,\, \,
  \Delta(F_i)=F_i\otimes K_i^{-1}+1\otimes F_i, \, \,
  \Delta(K_i)=K_i\otimes K_i,
\end{equation*}
\begin{equation*}
  S(E_i)=-K_i^{-1}E_i,\quad
  S(F_i)=-F_iK_i,\quad
  S(K_i)=K_i^{-1},
\end{equation*}
$$
\varepsilon(E_i)=\varepsilon(F_i)=0,\quad \varepsilon(K_i)=1.
$$

We need two important classes of $U_q \mathfrak{g}$-modules \cite{SSSV}. A
representation $\rho: U_q \mathfrak{g} \rightarrow {\rm End} V$ is called
{\it weight} (and $V$ is called a weight module, respectively), if $V$
admits a decomposition into a sum of weight subspaces
\begin{align*}
V = \bigoplus_{\boldsymbol \lambda} V_{\boldsymbol \lambda}, \qquad
V_{\boldsymbol \lambda}=\{v \in V |\text{ } \rho (K_j^{\pm 1})v = q_j^{\pm
\lambda_j}v, j=1,...,l \},
\end{align*}
where $\boldsymbol \lambda =(\lambda_1,...,\lambda_l) \in \mathbb Z^l$.
The subspace $V_{\lambda}$ is called a weight subspace of weight
$\lambda$.

Let $U_q \mathfrak{k} \subset U_q \mathfrak{g}$ be a Hopf subalgebra
generated by $E_i,F_i, i=1,...,l, i \neq l_0$ and $K^{\pm 1}_j, j=1,...,l$.
A finitely generated weight $U_q \mathfrak{g}$-module $V$ is called a
quantum Harish-Chandra module if $V$ is a sum of finite dimensional simple
$U_q \mathfrak{k}$-modules and $\dim \rm{Hom}_{U_q
\mathfrak{k}}(W,V)<\infty$ for every finite dimensional simple $U_q
\mathfrak{k}$-module $W$.

We restrict our consideration to quantum Harish-Chandra modules only.

 Equip $U_q \mathfrak{g}$ with a $*$-Hopf algebra structure
via the antilinear involution $*$:
$$
E_{l_0}^*=-K_{l_0}F_{l_0},\quad F_{l_0}^*=-E_{l_0}K_{l_0}^{-1},\quad
K_{l_0}^*=K_{l_0},
$$
$$
E_j^*=K_jF_j,\quad F_j^*=E_jK_j^{-1},\quad K_j^*=K_j,\quad j \not=l_0.
$$

Recall the notion of quantum analog of the algebra $\mathbb{C}[G]$ of
regular functions on $G$ \cite{Jantzen}.

Denote by $\mathbb{C}[G]_q \subset (U_q \mathfrak{g})^*$ the Hopf
subalgebra of all matrix coefficients of weight finite dimensional $U_q
\mathfrak{g}$-representations. $\mathbb{C}[G]_q$ is a $U_q
\mathfrak{g}$-module algebra:
$$
(\xi f)(\eta)=f(\eta \xi), \qquad \xi,\eta \in U_q \mathfrak{g}, f \in
\mathbb{C}[G]_q.
$$

The algebra $\mathbb{C}[G]_q$ is called the algebra of regular functions on
the quantum group $G$.

Introduce special notations for the elements of $\mathbb{C}[G]_q$
\cite{Soib}. Consider the finite dimensional simple $U_q
\mathfrak{g}$-module $L(\varpi_{l_0})$ with the highest weight
$\varpi_{l_0} \in P_+$. Equip it with an invariant scalar product
$(\cdot,\cdot)$ (as usual in the compact quantum group theory
\cite{Soib}). Following N. Reshetikhin and V. Lakshmibai (see
\cite{Resh}), choose nonzero vectors $\{v_{\mu}\} \in
L(\varpi_{l_0})_{\mu}$ for all weights $\mu \in W\varpi_{l_0}$. Let
$$
c_{\lambda,\mu}(\xi)= \frac{(\xi v_{\mu},v_{\lambda})}
{\|v_{\mu}\|\,\|v_{\lambda}\|} \qquad \mu,\lambda \in W\varpi_{l_0}.
$$

For brevity, put
\begin{equation*}
t= c_{\varpi_{l_0},-\varpi_{l_0}},\quad t'=c_{\varpi_{l_0},\varpi_{l_0}}.
\end{equation*}

Denote by $\mathbb{C}[X]_q \subset \mathbb{C}[G]_q$ the minimal $U_q
\mathfrak{g}$-module subalgebra generated by $t$.

Equip $\mathbb{C}[X]_q$ with an antilinear involution $*$ compatible with
the involution in $U_q \mathfrak{g}$, i.e.
\begin{equation*}
 (\xi f)^* = (S(\xi))^* f^*, \qquad
 \xi \in U_q\mathfrak{g},\; f \in \mathbb{C}[X]_q.
\end{equation*}

\begin{proposition}\label{homogen_*}
There exists a unique involution $*$ which equips the algebra
$\mathbb{C}[X]_q$ with the $(U_q \mathfrak{g}, *)$-module algebra
structure, such that $t^* = q^{(\varpi_{l_0},\rho)} t'$, where $\rho$ is
the half-sum of positive roots and the scalar product is defined by
$(\alpha_i,\alpha_j)=d_ia_{ij}$.
\end{proposition}

Let $x=tt^*$. Recall that $\mathbb{C}[X]_q$ is an integral domain.

\begin{proposition}
\begin{enumerate}
\item $x^{\mathbb{Z_+}}$ is an Ore set.
\item There exists a unique extension of the $(U_q \mathfrak{g},*)$-module
algebra structure from $\mathbb{C}[X]_q$ to the localization
$\mathbb{C}[X]_{q,x}$ of $\mathbb{C}[X]_q$ with respect to the
multiplicative set $x^{\mathbb{Z_+}}$.
\end{enumerate}
\end{proposition}

Equip $\mathbb{C}[X]_{q,x}$ with the $U_q \mathfrak{g}$-invariant $\mathbb
Z$-grading by $\deg t=1$. (E.g. $\deg ( c_{\varpi_{l_0},\mu}) =1$ for all
weights $\mu \in W\varpi_{l_0}$.)

Define the subalgebra $\mathbb{C}[S(\mathbb D)]_q =\{ f \in
\mathbb{C}[X]_{q,x} | \deg f=0\}$. This subalgebra inherits the $(U_q
\mathfrak{g},*)$-module algebra structure. This algebra is a quantum
analog of the algebra of regular functions on the Shilov boundary
$S(\mathbb D)$ of a bounded symmetric domain of tube type.

\section{\bf Examples of algebras of regular functions on the Shilov
boundaries}

We repeat below the general construction from the previous Sec. for the
special case $\mathfrak{g}=\mathfrak{sl}_N$. A Hermitian symmetric pair
$(\mathfrak{g},\mathfrak{k})$ is related to the irreducible bounded
symmetric domain of tube type only if $N=2n$ and $l_0=n$. The
corresponding bounded symmetric domain is the unit ball in the standardly
normed space of complex matrices $\mathbb D=\{A \in \mathrm{Mat}_{n,n}|
\,A A^* \leq I\}$. The Shilov boundary is isomorphic to the closed
$SU_{n,n}$-orbit of the Grassmanian $Gr_n(\mathbb C^{2n})$ under the Borel
embedding.

Recall the well known notation $U_q\mathfrak{su}_{n,n}$ for the $*$-Hopf
algebra $(U_q \mathfrak{sl}_{2n},*)$ and $\mathbb{C}[\mathrm{Mat}_{2n}]_q$
for the quantum $2n\times2n$-matrix space defined by the generators
$\{t_{ij}\}_{i,j=1,\ldots,2n}$ and the relations (cf. \cite{Drin_con})
\begin{align}\label{q_rel}
&t_{ik}t_{jk}=qt_{jk}t_{ik},\quad t_{ki}t_{kj}=qt_{kj}t_{ki},&i<j, \notag
\\ &t_{ij}t_{kl}=t_{kl}t_{ij},&i<k \; \& \; j>l,
\\ &t_{ij}t_{kl}-t_{kl}t_{ij}=(q-q^{-1})t_{ik}t_{jl},&i<k \; \& \; j<l.
\notag
\end{align}

$q$-minors are defined as follows:
\begin{equation}\label{q-min}
t_{IJ}^{\wedge k}\stackrel{\mathrm{def}}{=}\sum_{s \in
S_k}(-q)^{l(s)}t_{i_1j_{s(1)}}\cdots t_{i_kj_{s(k)}},
\end{equation}
for any $I=\{i_1,\ldots,i_k\},\;1 \le i_1<\ldots<i_k \le 2n,$
$J=\{j_1,\ldots,j_k\}, \;1 \le j_1<\ldots<j_k \le 2n$; here $l(s)$ denotes
the length of a permutation $s$. For brevity, put $\det_q
T=t_{\{1,...,2n\}\,\{1,...,2n\}}^{\wedge 2n}$. The algebra $\mathbb{C}[G]_q$
is obtained from $\mathbb{C}[\mathrm{Mat}_{2n}]_q$ by completing the list of
relations with the relation $\det_q T=1$.

The $U_q \mathfrak{sl}_{2n}$-module algebra $\mathbb{C}[X]_{q} \subset
\mathbb{C}[G]_{q}$ is generated by the $q$-minors $t^{\wedge
n}_{\{1,\ldots,n\}\,J}$ for all $J$ with $\operatorname{card} J=n$.

The compatible involution on $\mathbb{C}[X]_q$ is given by
$t^*=q^{\frac{n(n+1)}{2}}t'$. The elements $t\stackrel{\rm def}{=}
t^{\wedge n}_{\{1,\ldots,n\}\{n+1,\ldots,2n\}}$, $t'=t^{\wedge
n}_{\{n+1,\ldots,2n\}\{1,\ldots,n\}}$, and $x=tt^*$ quasi-commute with
$t_{ij}$ for all $i,j=1,\ldots,2n$.

Consider the $*$-algebra $\mathbb{C}[X]_{q,x}$. Obtain a description for
the subalgebra of zero-degree elements of $\mathbb{C}[X]_{q,x}$.

Let $\mathbb{C}[\mathrm{Mat}_n]_q \subset \mathbb{C}[X]_{q,x}$ be the
subalgebra generated by
\begin{align*}\label{emb}
z_a^b \stackrel{\mathrm{def}}{=} &t^{-1}t^{\wedge n}_{\{1,2,\ldots,n
\}J_{a \,b}}, \;\;\text{where} \quad J_{a \,b}=\{n+1,n+2,\ldots,2n
\}\setminus \{2n+1-b \}\cup \{a \}
\end{align*}
The defining relations between $z_a^b$ are similar to (\ref{q_rel})
\cite{V_qmbdic}:
\begin{align*}
&z_{a}^{b_1}z_{a}^{b_2}=qz_{a}^{b_2}z_{a}^{b_1}, &b_1<b_2,
\\&z_{a_1}^{b}z_{a_2}^{b}=qz_{a_2}^{b}z_{a_1}^{b}, &a_1<a_2,
\\&z_{a_1}^{b_1}z_{a_2}^{b_2}=z_{a_2}^{b_2}z_{a_1}^{b_1},
&b_1<b_2 \, \&\, a_1>a_2,
\\&z_{a_1}^{b_1}z_{a_2}^{b_2}-z_{a_2}^{b_2}z_{a_1}^{b_1}=
(q-q^{-1})z_{a_1}^{b_2}z_{a_2}^{b_1}, &b_1<b_2 \, \&\, a_1<a_2.
\end{align*}

Similarly to (\ref{q-min}), put $\mathbf z ^{\wedge k
\{b_1,\ldots,b_k\}}_{\,\,\,\,\,\,\,\{a_1,\ldots,a_k\}} \stackrel {\rm
def}{=} \sum_{s \in S_k} (-q)^{l(s)} z^{b_{s(1)}}_{a_1}\ldots
z^{b_{s(k)}}_{a_k}$, where $a_1 < \ldots < a_k,\,\, b_1<\ldots<b_k$ and
$\det_q \mathbf z\stackrel{\mathrm{def}}{=}\mathbf
{z}_{\,\,\,\,\,\,\,\{1,\ldots,n\}}^{\wedge n\{1,\ldots,n\}}$. It is easy
to prove that $\det_q \mathbf z=t^{-1}t'$. The algebra
$\mathbb{C}[S(\mathbb D)]_q$ is isomorphic to the localization of the
algebra $\mathbb{C}[\mathrm{Mat}_n]_q$ with respect to the Ore system
$(\det_q \mathbf z)^{\mathbb{Z}_+}$ and
\begin{equation*}
(z_a^b)^*=(-q)^{a+b-2n}({\rm det}_q \mathbf z)^{-1} {\rm det}_q \mathbf
z_a^b ,
\end{equation*}
where $\det_q \mathbf z_a^b$ is the $q$-determinant of the matrix derived
from $\mathbf z$ by deleting the line $b$ and the column $a$.

The $U_q \mathfrak{su}_{n,n}$-action on $\mathbb{C}[S(\mathbb D)]_q$ can
be described explicitly:
\begin{align*}
K_n^{\pm 1}z_a^b=&
\begin{cases}
q^{\pm 2}z_a^b,&a=n \;\&\;b=n
\\ q^{\pm 1}z_a^b,&a=n \;\&\;b \ne n \;\text{or}\;
a \ne n \;\&\; b=n
\\ z_a^b,&\text{\rm otherwise},
\end{cases}
\\ F_nz_a^b=q^{1/2}\cdot
\begin{cases}
1,& a=n \;\& \;b=n
\\ 0,&\text{\rm otherwise},
\end{cases}\quad
& \quad E_nz_a^b=-q^{1/2}\cdot
\begin{cases}
q^{-1}z_a^nz_n^b,&a \ne n \;\&\;b \ne n
\\ (z_n^n)^2,& a=n \;\&\;b=n
\\ z_n^nz_a^{b},&\text{\rm otherwise}
\end{cases}
\end{align*}
and for all $k \ne n$
\begin{align*}
\!\!\!\!\!\!K_k^{\pm 1}z_a^b=
\begin{cases}
\!q^{\pm 1}z_a^b,\; k<n \;\&\;a=k \;\text{or} \;k>n \;\&\;b=2n-k,
\\ \!q^{\mp 1}z_a^b,\; k<n \;\&\;a=k+1 \;\text{or} \;k>n \;\&\;b=2n-k+1,
\\ \!z_a^b,\; \text{\rm otherwise},
\end{cases}
\\ \!\!\!F_kz_a^b=\!q^{1/2}\!\cdot
\begin{cases}
\!z_{a+1}^b, \;k<n \;\&\;a=k,
\\ \!z_a^{b+1}, \;k>n \;\&\;b=2n-k,
\\ \!0, \;\; \text{\rm otherwise},
\end{cases}
\!\!\! E_kz_a^b=\!q^{-1/2}\!\cdot
\begin{cases}
\!z_{a-1}^b, \;k<n \;\&\; a=k+1,
\\ \!z_a^{b-1}, \;k>n \;\&\;b=2n-k+1,
\\ \!0, \;\;\;\;\text{\rm otherwise}.
\end{cases}
\end{align*}

The $*$-algebra $\mathbb{C}[S(\mathbb D)]_q$ is isomorphic to the algebra
of regular functions on the quantum group $U_n$ \cite{FRT}. The latter
algebra was equipped with the $U_q \mathfrak{k}$-module algebra structure.
Our general construction gives us hidden symmetry (additional structure of
$U_q \mathfrak{g}$-module algebra).

\medskip

Now turn to type $C_n$ Lie algebras. Fix $\mathfrak{g}=\mathfrak{sp}_{2n}$.
The pair $(\mathfrak{g},\mathfrak{k})$ is Hermitian symmetric only if
$\mathfrak{k}=\mathfrak{gl}_n$. In this case the bounded symmetric domain is
the unit ball in the standardly normed space of complex symmetric $n \times
n$-matrices. The $*$-Hopf algebra $(U_q \mathfrak{sp}_{2n},*)$ is a quantum
analog of the universal enveloping algebra $U
\mathfrak{sp}_{2n}(\mathbb{R})$.

We need an algebra $\mathbb{C}[\mathrm{Mat}^{sym}_n]_q$ defined by the
generators $z_{ij}$ for $1 \leq j \leq i \leq n$ and the
relations\footnote{This is a quantum analog of the polynomial algebra on
complex symmetric $n \times n$-matrices.}
\begin{align*}
 z_{ii}z_{ki}&=q^2z_{ki}z_{ii}, &i < k
\\ z_{ki}z_{kk}&=q^2z_{kk}z_{ki}, & i<k
\\ z_{ij}z_{ik}&=qz_{ik}z_{ij}, &j<k<i
\\ z_{ij}z_{kj}&=q z_{kj}z_{ij}, &j<i<k
\\ z_{ij}z_{kl}&=z_{kl}z_{ij}, &j<l \leq k<i
\\ z_{ii}z_{jj}&=z_{jj}z_{ii}+q(q^2-q^{-2})z^2_{ji}, &i<j
\\ z_{ii}z_{jk}&=z_{jk}z_{ii}+(q^2-q^{-2})z_{ki}z_{ji}, &i<k<j
\\ z_{ik}z_{jj}&=z_{jj}z_{ik}+(q^2-q^{-2})z_{jk}z_{ji}, &k<i<j
\\ z_{ij}z_{kl}&=z_{kl}z_{ij}+(q-q^{-1})(qz_{li}z_{kj}+z_{ki}z_{lj}), &j<i<l<k
\\ z_{ij}z_{kl}&=z_{kl}z_{ij}+(q-q^{-1})z_{il}z_{kj}, &j<l<i<k
\\ z_{ij}z_{kl}&=qz_{kl}z_{ij}+(q-q^{-1})z_{il}z_{kj}, &j<i=l<k
\end{align*}
{\it Remark.} The algebra generated by linear functionals on $\{z_{ij}\}$
was introduced in \cite{Kam}, where it was considered as a $U_q
\mathfrak{k}$-module algebra.

\medskip

Denote a quantum analog of determinant of symmetric matrix
\begin{equation*}
{\rm det}^{sym}_q \mathbf z\stackrel {\rm def}{=} \sum_{s \in S_n}
(-q)^{-l(s)} q^{n-\sum_{i=1}^n\delta_{is(i)}}z_{s(n)\,n}\cdot \ldots \cdot
z_{s(1)\,1}
\end{equation*}
with $z_{kl}=q^{-2}z_{lk}$ for $l> k$.

Similarly to the case $A_n$, the algebra $\mathbb{C}[S(\mathbb D)]_q$ is
isomorphic to the localization of the algebra
$\mathbb{C}[\mathrm{Mat^{sym}}_n]_q$ with respect to the Ore system
$(\det^{sym}_q \mathbf z)^{\mathbb Z_+}$ with the involution
\begin{equation*}
z_{nn}^*= {\rm det}^{sym}_q \mathbf z_{nn} ({\rm det}^{sym}_q \mathbf
z)^{-1},
\end{equation*}
where $\det^{sym}_q \mathbf z_{nn}$ is the $q$-determinant of the matrix
derived from $\mathbf z$ by deleting the line $n$ and the column $n$.

From the general construction, $\mathbb{C}[S(\mathbb D)]_q$ is equipped
with the $(U_q \mathfrak{sp}_{2n},*)$-module algebra structure. The $U_q
\mathfrak{sp}_{2n}$-action is described explicitly via the following:
\begin{equation*}
K_n z_{ij}=
\begin{cases}
q^4 z_{ij}, & i=j=n,
\\q^2 z_{ij}, & i=n>j \; \text{or}\; j=n>i,
\\ z_{ij}, & \text{otherwise}.
\end{cases}
\end{equation*}
\begin{equation*}
F_n z_{ij}= q \begin{cases} z_{ij}, & i=j=n,
\\ 0, & \text{otherwise}.
\end{cases}
\qquad E_n z_{ij}= -q
\begin{cases}
z_{nn}z_{ij}, & i=n \geq j \; \text{or} \;j=n \geq i,
\\ q^{-1}z_{ni}z_{nj}, & \text{otherwise}.
\end{cases}
\end{equation*}

and for $k \neq n$ and $i \geq j$
\begin{equation*}
E_k z_{ij}= q^{-1/2}
\begin{cases}
 (q+q^{-1})z_{i \, j-1} & i=j=k+1,
\\ z_{i-1 \, j} & i=k+1 \,\&\, j<k+1,
\\ z_{i \, j-1} & i>k+1 \,\&\, j=k+1,
\\ 0 & \text{otherwise}.
\end{cases}
\end{equation*}
\begin{equation*}
 F_k z_{ij} = q^{1/2} \begin{cases}
(q+q^{-1})z_{i+1\,j} \,&\, i=j=k,
\\ z_{i+1\,j} & i=k \,\&\, j<k,
\\ z_{i\,j+1} & i>k \&\, j=k,
\\ 0 & \text{otherwise}.
\end{cases}
\end{equation*}
\begin{equation*}
K_k z_{ij} =\begin{cases} q^2 z_{ij}& i=j=k,
\\ q^{-2}z_{ij} & i=j=k+1,
\\ qz_{ij} & (i=k \,\&\, j<k) \,\text{or}\, (i>k+1 \,\&\, j=k),
\\ q^{-1}z_{ij} & (i=k+1 \,\&\, j<k) \,\text{or}\, (i>k+1 \,\& \,
j=k+1),
\\ z_{ij} & \text{otherwise}.
\end{cases}
\end{equation*}

Here is another interesting fact about the Shilov boundaries.

In the cases $A_n$ and $C_n$ there are points on the Shilov boundaries,
i.e. $(U_q \mathfrak{g},*)$-morphisms $p:\mathbb{C}[S(\mathbb D)]_q
\rightarrow \mathbb C$. The respective $*$-morphisms can be rebuilt from
the formulas
\begin{equation*}
p(z_a^b)= \begin{cases} q^{n-a}, & a=b,
\\ 0, &a \neq b, \end{cases}
\quad p(z_{ij})= \begin{cases} q^{n-i}, & i=j,
\\ 0, & i \neq j. \end{cases}
\end{equation*}
Using these points, we can construct in both cases a $(U_q
\mathfrak{g},*)$-morphism $i:\mathbb{C}[S(\mathbb D)]_q \rightarrow
\mathbb{C}[K]_q$ such that the following diagram is commutative
$$\xymatrix{\mathbb{C}[S(\mathbb D)]_q \ar[dr]^p \ar[r]^i & \mathbb{C}[K]_q
\ar[d]^{\varepsilon}\\ & \mathbb C}$$

The subalgebra $i(\mathbb{C}[S(\mathbb D)]_q)$ in these cases admits a
description in the spirit of M. Noumi's paper \cite{Stok}.

\section{Degenerate principal series of quantum Harish-Chandra
modules related to the Shilov boundary for the case $A_n$} \label{probl}

In this section we investigate a quantum analog of the degenerate
principal series of $U_q\mathfrak{su}_{n,n}$-modules related to the Shilov
boundary of the quantum $n \times n$-matrix unit ball. We give necessary
and sufficient conditions for the representations to be irreducible and
unitarizable.

In this section we provide $q$-analogs of classical results obtained by
K.D. Johnson, S. Sahi, G. Zhang, R.E. Howe, and E.-C. Tan
\cite{Howe,Jons1,Jons2,Sahi,Zhang}. Another degenerate principal series is
considered in the A. Klimyk and S. Pakuliak paper \cite{KL_Pak}. Our
results are quantum analogs of results from \cite{Lee}. More detailed
results and proofs are given in \cite{Deg_ser}.

Assume first that $\alpha,\beta \in \mathbb Z$. Define a representation
$\pi_{\alpha,\beta}:U_q\mathfrak{sl}_{2n} \rightarrow {\rm End} (\mathbb
C[S(\mathbb{D})]_q)$ as follows:
\begin{equation*}
\pi_{\alpha,\beta}(\xi)f=(\xi \cdot (f (t')^{\alpha} t^{\beta}))
t^{-\beta} (t')^{-\alpha} = (\xi \cdot (f ({\rm det}_q \mathbf z)^{\alpha}
t^{\beta+\alpha})) t^{-\alpha-\beta} ({\rm det}_q \mathbf z)^{-\alpha}
\end{equation*}
for every $\xi \in U_q\mathfrak{sl}_{2n},f \in \mathbb
C[S(\mathbb{D})]_q$. For each $\lambda \in \mathbb Z$ we have
$$
E_j t^{\lambda}=0,\,\, F_jt^{\lambda}=0,\,\, K_jt^{\lambda}=t^{\lambda},
\quad j=1,\ldots,2n-1, \; j\ne n
$$
$$
E_nt^\lambda=q^{-3/2}\frac{1-q^{-2\lambda}}{1-q^{-2}}z^n_n t^\lambda,
\quad F_nt^\lambda=0, \quad K_n^{\pm 1}t^\lambda=q^{\mp \lambda}t^\lambda,
$$
$$
E_j ({\rm det}_q \mathbf z)^\lambda=0,\,\, F_j ({\rm det}_q \mathbf
z)^\lambda=0,\,\, K_j ({\rm det}_q \mathbf z)^\lambda=({\rm det}_q \mathbf
z)^\lambda, \quad j=1,\ldots,2n-1, \; j\ne n
$$
$$
K_n^{\pm 1}(({\rm det}_q \mathbf z)^\lambda)=q^{\pm 2 \lambda}({\rm det}_q
\mathbf z)^\lambda, \quad E_n(({\rm det}_q \mathbf z)^\lambda) = -q^{1/2}
\frac{1-q^{2\lambda}}{1-q^2} z^n_n ({\rm det}_q \mathbf z)^\lambda,
$$
$$
F_n(({\rm det}_q \mathbf z)^\lambda) =
q^{1/2}\frac{1-q^{-2\lambda}}{1-q^{-2}} \mathbf z^{\wedge
n-1}_{\{1,\ldots,n-1\}\{1,\ldots,n-1\}} ({\rm det}_q \mathbf
z)^{\lambda-1},\quad \lambda \ne 0.
$$

From these equalities we see that for each $\xi \in
U_q\mathfrak{sl}_{2n}$, $f \in \mathbb C[S(\mathbb{D})]_q$ the vector
function $p_{f,\xi}(q^{\alpha},q^{\beta}) \stackrel{\rm
def}=\pi_{\alpha,\beta}(\xi)(f)$ is a Laurent polynomial of the variables
$q^{\alpha},\,q^{\beta}$. These Laurent polynomials are uniquely defined
by their values on the set $\{(q^{\alpha},q^{\beta})| \,\alpha,\,\beta \in
\mathbb{Z}\}$ and deliver the canonical "analytic continuation" for
$\pi_{\alpha,\beta}(\xi)(f)$ to $(\alpha,\beta) \in \mathbb{C}^{2}$.

Let $(\alpha,\beta) \in \mathbb{C}^{2}$. Define a representation
$\pi_{\alpha,\beta}(\xi)(f) \stackrel{\rm def}=p_{f,\xi}(\alpha,\beta)$.
To prove that the representation $\pi_{\alpha,\beta}$ is well defined for
$(q^{\alpha},\,q^{\beta}) \in \mathbb{C}^{2}$, it is sufficient to verify
some identities for Laurent polynomials. These identities hold for
$\alpha,\,\beta \in \mathbb{Z}$.

Introduce a ''deformation parameter'' $h$ by the equality $q=e^{-h/2}$.
Clearly, if $\alpha_1=\alpha_2+i\frac{2\pi}{h}$ and
$\beta_1=\beta_2+i\frac{2\pi}{h}$, then
$\pi_{\alpha_1,\beta_1}=\pi_{\alpha_2,\beta_2}$. Then it suffices to
consider $(\alpha,\,\beta) \in D$, where
\begin{equation*}
D=\{ (\alpha,\beta) \in \mathbb C^2 \, |\, 0 \leq \operatorname{Im} \alpha
< \frac{2\pi}{h}, \, 0 \leq \operatorname{Im} \beta < \frac{2\pi}{h}\}.
\end{equation*}

It is clear that $\pi_{\alpha,\beta}$ defines a Harish-Chandra module if and
only if $q^{\alpha-\beta} \in q^{\mathbb Z}$. Note that for any complex
$\alpha,\beta$ such that $0 \leq \operatorname{Im} \alpha < \frac{2\pi}{h},
\, 0 \leq \operatorname{Im} \beta < \frac{2\pi}{h}$, the statements
$q^{\alpha-\beta} \in q^{\mathbb Z}$ and $\alpha-\beta \in \mathbb Z$ are
equivalent.

\subsection {Equivalence of the representations}\label{equi}

In this subsection we obtain a parameter set $\mathcal{D}$ and prove that
each representation of the degenerate principal series is equivalent to a
representation $\pi_{\alpha,\beta}$ for some $(\alpha,\beta) \in
\mathcal{D}$.

The representations $\pi_{\alpha,\beta}$ and $\pi_{\alpha-1,\beta+1}$ are
equivalent for all $\alpha,\beta$. The corresponding intertwining operator
$T: \mathbb C[S(\mathbb{D})]_q \rightarrow \mathbb C[S(\mathbb{D})]_q$ is
defined as follows: $T(f)=f \cdot (\det_q \mathbf z)^{-1}$ for every $f
\in \mathbb C[S(\mathbb{D})]_q$. Indeed, since for each $f \in \mathbb
C[S(\mathbb{D})]_q,\,\xi \in U_q \mathfrak{sl}_{2n}$
\begin{multline*}
\pi_{\alpha-1,\beta+1}(\xi)(f)=(\xi \cdot (f ({\rm det}_q \mathbf
z)^{\alpha-1} t^{\beta+\alpha})) t^{-\alpha-\beta} ({\rm det}_q \mathbf
z)^{1-\alpha}=
\\ (\xi \cdot (f ({\rm det}_q \mathbf z)^{-1} ({\rm det}_q \mathbf z)^{\alpha}
t^{\beta+\alpha})) t^{-\alpha-\beta} ({\rm det}_q \mathbf z)^{-\alpha}
({\rm det}_q \mathbf z)= \pi_{\alpha,\beta}(\xi)(f({\rm det}_q \mathbf
z)^{-1}){\rm det}_q \mathbf z,
\end{multline*}
we see that $T$ intertwines the representations $\pi_{\alpha,\beta}$ and
$\pi_{\alpha-1,\beta+1}$. Therefore without loss of generality we can
assume that $\alpha,\beta \in \mathcal{D}$, where
\begin{equation*}
\mathcal{D}=\{ (\alpha,\beta) \in \mathbb
C^2 \, |\, \alpha-\beta \in \{0,1\},\, 0 \leq \operatorname{Im} \alpha <
\frac{2\pi}{h}, 0 \leq \operatorname{Im} \beta < \frac{2\pi}{h}\}.
\end{equation*}

\begin{proposition} \label{pr_eqv}
If $\alpha,\beta \not \in \mathbb Z$, then the representations
$\pi_{\alpha,\beta}$
 and $\pi_{-n-\beta,-n-\alpha}$ are equivalent.
\end{proposition}
The proof follows from the explicit formulas for the intertwining
operators given in Section \ref{eqv}.

If $\alpha,\beta \in \mathbb Z$, then the representations
$\pi_{\alpha,\beta}$ and $\pi_{-n-\beta,-n-\alpha}$ are not equivalent.
This fact follows from the statement that only one of the representations
$\pi_{\alpha,\beta}$ and $\pi_{-n-\beta,-n-\alpha}$ has a finite
dimensional subrepresentation for integral $\alpha,\beta$. An explanation
of this fact is given in the end of Section \ref{reduce}.

Introduce an equivalence relation on $\mathcal{D}$. The equivalence class
of $(\alpha,\beta)$ consists of one point for $\alpha,\beta \in \mathbb Z$
and two points for $\alpha,\beta \not \in \mathbb Z$:
\begin{equation*}
(\alpha_1,\beta_1) \sim (\alpha_2,\beta_2), \quad \text{iff} \quad
\begin{cases} \alpha_1=-n-\beta_2,\,\beta_1=-n-\alpha_2 \,\,\text{for} \,\,
\operatorname{Im} \alpha_1 = \operatorname{Im} \alpha_2 =0,
\\ \alpha_1=\frac{2 \pi i}{h}-n-\beta_2,\,\beta_1=\frac{2 \pi i}{h}-n-\alpha_2,
\,\, \text{otherwise}.
\end{cases}
\end{equation*}
\begin{proposition}\label{clas_eqv}
The set of equivalence classes $\mathcal{D} \diagup \sim$ is in the
one-to-one correspondence $(\alpha,\beta) \mapsto \pi_{\alpha,\beta}$ with
the set of equivalence classes of representations of the degenerate
principal series.
\end{proposition}
{\bf Proof.} By the above, each representation of the degenerate principal
series is equivalent to the representation $\pi_{\alpha,\beta}$ for some
$(\alpha,\beta) \in \mathcal{D}$.

Prove that the representations $\pi_{\alpha_1,\beta_1}$ and
$\pi_{\alpha_2,\beta_2}$, with $(\alpha_1,\beta_1),(\alpha_2,\beta_2) \in
\mathcal{D}$, are equivalent if and only if $(\alpha_1,\beta_1) \sim
(\alpha_2,\beta_2)$. For that, we calculate the action of a central
element $C \in U_q\mathfrak{sl}_{2n}^{ext}$ (see \cite{Klimyk} for the
definition). It can be proved that $\pi_{\alpha,\beta}(C)$ is a scalar
operator for all $\alpha,\beta \in \mathcal{D}$.

From \cite{Drin} it follows that there exists a unique central element $C$
which acts on the $U_q\mathfrak{sl}_{2n}$-highest vector $v^{high}$ with
weight $\lambda$ 
as follows:
\begin{equation*}
C({v^{high}})=\sum \limits_{j=0}^{2n-1}
q^{-2(\mu_j,\lambda+\rho)}v^{high},
\end{equation*}
where $\mu_0=\varpi_1, \;\mu_j=-\varpi_j+\varpi_{j+1}$ for
$j=1,\ldots,2n-2$, $\mu_{2n-1}=-\varpi_{2n-1}$.

First let $\alpha,\beta$ be integers. It can be proved that
\begin{equation*}
\pi_{\alpha,\beta}(C)({\rm det}_q \mathbf z)^{\beta}=
4\operatorname{ch}\frac{h}{2} (\alpha+\beta+n) (\sum_{j=0}^{n-1}
\operatorname{ch} \frac{h}{2}j) ({\rm det}_q \mathbf z)^{\beta}.
\end{equation*}
Hence $\pi_{\alpha,\beta}(C)=4\operatorname{ch}\frac{h}{2}
(\alpha+\beta+n) (\sum_{j=0}^{n-1} \operatorname{ch} \frac{h}{2}j)\cdot
\mathrm{Id}$ for all $(\alpha,\beta) \in \mathcal{D}$.

Suppose that $\pi_{\alpha_1,\beta_1}$ and $\pi_{\alpha_2,\beta_2}$ are
equivalent. Equivalent representations have the same weight lattice.
Therefore $(\alpha_1-\beta_1)-(\alpha_2-\beta_2) \in 2\mathbb Z$. Since
$(\alpha_1,\beta_1),(\alpha_2,\beta_2) \in \mathcal{D}$, we see that
$(\alpha_1-\beta_1)-(\alpha_2-\beta_2)=0$.

Then the equivalent representations $\pi_{\alpha_1,\beta_1}$ and
$\pi_{\alpha_2,\beta_2}$ have the same values of central characters, which
means that
$$
(\operatorname{ch}\frac{h}{2}(\alpha_1+\beta_1+n)-
\operatorname{ch}\frac{h}{2}(\alpha_2+\beta_2+n)) \sum_{j=0}^{n-1}
\operatorname{ch} \frac{h}{2}j=0
$$
Since $0 \leq \operatorname{Im} \alpha_1 < \frac{2\pi}{h}$, $0 \leq
\operatorname{Im} \beta_1 < \frac{2\pi}{h}$, $0 \leq \operatorname{Im}
\alpha_2 < \frac{2\pi}{h}$, $0 \leq \operatorname{Im} \beta_2 <
\frac{2\pi}{h}$, we have that $\alpha_1+\beta_1=\alpha_2+\beta_2$, or
$\alpha_1+\beta_1=-\alpha_2-\beta_2-2n$, or
$\alpha_1+\beta_1=-\alpha_2-\beta_2-2n-\frac{4\pi i}{h}$. If
$\alpha_1+\beta_1=\alpha_2+\beta_2$, then $\alpha_1=\alpha_2$ and
$\beta_1= \beta_2$. For any fixed non-integral $\alpha_1,\beta_1$ there is
a unique pair $(\alpha_2,\beta_2) \in \mathcal{D}$ such that either
$\alpha_1+\beta_1=-\alpha_2-\beta_2-2n$ or
$\alpha_1+\beta_1=-\alpha_2-\beta_2-2n-\frac{4\pi i}{h}$, and
$(\alpha_1,\beta_1) \sim (\alpha_2,\beta_2)$. Although for integral
parameters $\pi_{\alpha_1,\beta_1}$ and $\pi_{\alpha_2,\beta_2}$ are not
equivalent, because the only one of them has a finite-dimensional
subrepresentation. This can be deduced from Corollary \ref{cor_4}. Thus
each equivalence class in $\mathcal{D}$ is assigned to a unique
equivalence class of representations of the degenerate principal series
$\pi_{\alpha,\beta}$.\hfill $\square$

\subsection{Reducibility of $\pi_{\alpha,\beta}$}\label{reduce}

Let $U_q \mathfrak{k}_{ss}\subset U_q \mathfrak{sl}_{2n}$ be the Hopf
subalgebra generated by $E_j, F_j, K_j^{\pm 1}$, $j=1,\ldots,2n-1$, $j
\neq n$ and $U_q \mathfrak{k}\subset U_q \mathfrak{sl}_{2n}$ be the Hopf
subalgebra generated by $K_n^{\pm 1}$ and $U_q \mathfrak{k}_{ss}$.

Note that $\pi_{\alpha,\beta}|_{U_q \mathfrak{k}_{ss}}$ does not depend on
$\alpha,\beta$. The following preliminary result on reducibility of
$\pi_{\alpha,\beta}$ is well known in the classical case. For brevity,
set\footnote{Note that, obviously, $\mathbf z^{\wedge n}=\det_q \mathbf
z$.}
$$
\mathbf z^{\wedge k} = \mathbf z^{\wedge k
\{1,\ldots,k\}}_{\,\,\,\,\,\,\{1,\ldots,k\}}.
$$

Introduce the following notation: $\widehat{K}=\{\overline{\mathbf
k}=(k_1,\ldots,k_n) \in \mathbb Z^n|\,k_1 \geq k_2 \geq \ldots \geq
k_n\}$, $\mathbf{e_j}=(0,\ldots,\overset{j}{1},\ldots,0) \in \mathbb Z^n$.

\begin{proposition} \label{K_reduce}
The representation space $\mathbb C[S(\mathbb D)]_q$ for
$\pi_{\alpha,\beta}$ splits into a sum of simple pairwise non-isomorphic
$U_q \mathfrak{k}$-modules as follows:\footnote{These isotypic components
are $U_q \mathfrak{k}_{ss}$-isomorphic. However, they are not $U_q
\mathfrak{k}$-isomorphic, since the action of $\pi_{\alpha,\beta}(K_n)$
depends on $\alpha,\beta$.}
\begin{equation*}
\mathbb C[S(\mathbb D)]_q=\bigoplus_{\overline{\mathbf k} \in \widehat{K}}
V_{\overline{\mathbf k}}, \quad \text{with} \quad V_{\overline{\mathbf
k}}=\pi_{\alpha,\beta}(U_q\mathfrak{k}) \cdot v^h_{\overline{\mathbf k}}
\quad \text{and} \quad v^h_{\overline{\mathbf k}}=(\mathbf z^{\wedge
1})^{k_1-k_2}\ldots(\mathbf z^{\wedge n-1})^{k_{n-1}-k_n}(\mathbf
z^{\wedge n})^{k_n}.
\end{equation*}
\end{proposition}

\medskip
{\it Remark.} It can be easily verified that $v^h_{\overline{\mathbf k}}$
is a $U_q\mathfrak{k}$-highest vector with weight
$(k_1-k_2,\ldots,k_{n-1}-k_n,2k_n+\alpha-\beta,k_{n-1}-k_n,\ldots,k_1-k_2)$.
Then the highest weight of the simple $U_q\mathfrak{k}$-module
$V_{\overline{\mathbf k}}$ is equal to
$(k_1-k_2,\ldots,k_{n-1}-k_n,2k_n+\alpha-\beta,k_{n-1}-k_n,\ldots,k_1-k_2)$.

\begin{proposition} \label{pr_rez}
The representation $\pi_{\alpha,\beta}$ is irreducible if and only if
$\alpha,\beta$ satisfy the following equivalent conditions:\footnote{Since
$\alpha-\beta \in \mathbb Z$, these conditions are equivalent.}
\begin{equation*}
1.\,\,\alpha \not \in \mathbb Z; \qquad 2.\,\,\beta \not \in \mathbb Z.
\end{equation*}
\end{proposition}

\medskip

Now suppose $\alpha,\,\beta \in \mathbb Z$. We investigate reducibility
and proper subrepresentations of $\pi_{\alpha,\,\beta}$. We describe
results with figures as in \cite{Howe,Lee}.

Each $U_q\mathfrak{k}$-isotypic component $V_{\overline{\mathbf k}}$ is
assigned to the point $(k_1,\ldots,k_n) \in \mathbb R^n$. Thus
$\widehat{K}$ is assigned to the set $\mathbf
K^{+}=\{(k_1,\ldots,k_n)\,|\,k_1 \geq \ldots \geq k_n \} \subset \mathbb
R^n$. Consider $2n$ hyperplanes:
$$\mathcal{L}_j^{+}:k_j=\beta+j-1; \qquad \mathcal{L}_j^{-}: k_j=-\alpha-n+j.$$

These hyperplanes are parallel to the coordinate axis and pass through
points with integral coordinates. The distance between $\mathcal{L}_j^{+}$
and $\mathcal{L}_j^{-}$ is equal to $\alpha+\beta+n-1$.

\medskip

Investigate the example $n=2$. In this case $\mathcal{L}_j^{\pm},\,j=1,2$,
are just lines on the plane $\mathbb{R}^2$, parallel to coordinate axis.
Consider different values of $\alpha+\beta$.

Case 1. $\alpha+\beta \geq 0$. In this case $\mathcal{L}_1^{+}$ lies to
the right of $\mathcal{L}_1^{-}$, $\mathcal{L}_2^{+}$ lies higher than
$\mathcal{L}_2^{-}$ (see Fig.1). The intersection point of
$\mathcal{L}_1^{+}$ and $\mathcal{L}_2^{-}$ has the coordinates
$(\beta,\,-\alpha)$ and belongs to $\mathbf K^+$. Arrows attached to
$\mathcal{L}_j^{\pm}$ show the direction where $\pi_{\alpha,\beta}$
"moves" the isotypic components. There exists a unique simple submodule
$V^s=\bigoplus \limits_{\{\overline{\mathbf k} \in \widehat{K}|k_1 \leq
\beta,\,k_2 \geq -\alpha\}} V_{\overline{\mathbf k}}$ in $\mathbb
C[S(\mathbb D)]_q$.

\medskip

\begin{picture} (300,330)(0,0)
\put(15,20){\dottedline{2}(0,0)(0,300)} \put(15,325){\vector(0,1){5}}
\put(5,30){\dottedline{2}(0,0)(300,0)} \put(310,30){\vector(1,0){5}}
\put(5,100){\line(1,0){310}} \put(5,180){\line(1,0){310}}
\put(65,20){\line(0,1){310}} \put(125,20){\line(0,1){310}}
\put(65,40){\vector(1,0){15}} \put(300,180){\vector(0,-1){15}}
\put(300,100){\vector(0,1){15}} \put(125,40){\vector(-1,0){15}}
\put(10,25){\line(1,1){295}} \put(5,305){$k_2$} \put(305,20){$k_1$}
\put(70,20){$-\alpha-1$} \put(130,20){$\beta$} \put(0,105){$-\alpha$}
\put(0,185){$\beta+1$}

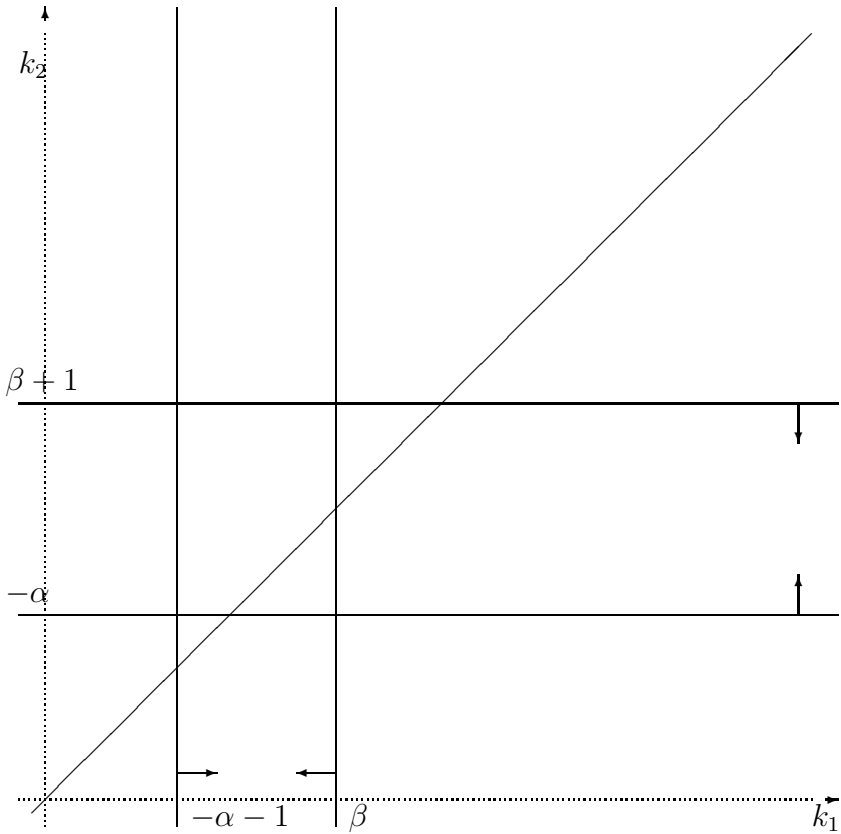
\put(100,0){Fig.1. The structure of $\pi_{\alpha,\beta}$ with
$\alpha+\beta \geq 0$.}
\end{picture}

\medskip

Case 2. $\alpha+\beta=-1$. In this case $\mathcal{L}_1^{+}$ and
$\mathcal{L}_1^{-}$, $\mathcal{L}_2^{+}$ and $\mathcal{L}_2^{-}$ coincide.
The intersection point of $\mathcal{L}_1^{+}$ and $\mathcal{L}_2^{+}$ does
not belong to $\mathbf K^+$ (Fig.2). There are two simple submodules in
$\mathbb C[S(\mathbb D)]_q$: $V^s_1= \bigoplus
\limits_{\{\overline{\mathbf k} \in \widehat{K}| k_1=
-1-\alpha\}}V_{\overline{\mathbf k}}$ and $V^s_2= \bigoplus \limits_{
\{\overline{\mathbf k} \in \widehat{K}| k_2 =-\alpha\}}
V_{\overline{\mathbf k}}$.

\begin{picture} (300,330)(0,0)
\put(15,20){\dottedline{2}(0,0)(0,300)} \put(15,315){\vector(0,1){5}}
\put(5,30){\dottedline{2}(0,0)(300,0)} \put(310,30){\vector(1,0){5}}
\put(5,130){\line(1,0){310}} \put(105,10){\line(0,1){310}}
\put(90,40){\vector(1,0){15}} \put(300,145){\vector(0,-1){15}}
\put(300,115){\vector(0,1){15}} \put(120,40){\vector(-1,0){15}}
\put(10,25){\line(1,1){295}} \put(5,310){$k_2$} \put(305,20){$k_1$}
\put(110,20){$\beta$} \put(0,135){$-\alpha$}

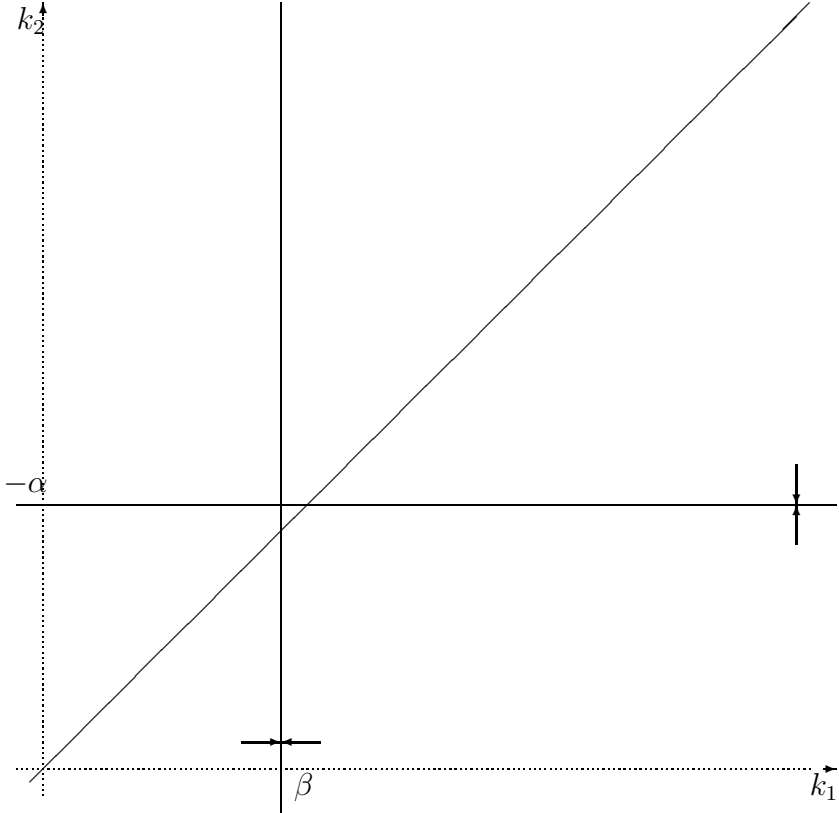
\put(100,0){Fig.2. The structure of $\pi_{\alpha,\beta}$ with
$\alpha+\beta=-1$.}
\end{picture}

\medskip

Case 3. $\alpha+\beta=-2$. In this case $\mathcal{L}_1^{+}$ lies to the
left of $\mathcal{L}_1^{-}$, $\mathcal{L}_2^{+}$ lies lower than
$\mathcal{L}_2^{-}$. However, $\mathcal{L}_1^{-}$ and $\mathcal{L}_2^{+}$
meet at the point with coordinates $(-\alpha-1,\,\beta+1)$ (see Fig.3).
Besides, the distance between $\mathcal{L}_j^{+}$ and
$\mathcal{L}_j^{-}$is 1. This shows that $\mathbb C[S(\mathbb D)]_q$ is a
direct sum of three submodules:
$$
V^s_1= \bigoplus \limits_{\{\overline{\mathbf k} \in \widehat{K}|k_1 \leq
\beta\}} V_{\overline{\mathbf k}}, \quad V^s_2= \bigoplus
\limits_{\{\overline{\mathbf k} \in \widehat{K}|k_2 \geq -\alpha\}}
V_{\overline{\mathbf k}},\quad V^s_3= \bigoplus
\limits_{\{\overline{\mathbf k} \in \widehat{K}|k_1 \geq -\alpha-1,k_2
\leq \beta+1\}} V_{\overline{\mathbf k}}.
$$

\begin{picture}(300,300)(0,0)
\put(15,20){\dottedline{2}(0,0)(0,300)} \put(15,325){\vector(0,1){5}}
\put(5,30){\dottedline{2}(0,0)(300,0)} \put(310,30){\vector(1,0){5}}
\put(5,120){\line(1,0){310}} \put(5,200){\line(1,0){310}}
\put(65,20){\line(0,1){310}} \put(145,20){\line(0,1){310}}
\put(145,40){\vector(1,0){15}} \put(285,120){\vector(0,-1){15}}
\put(285,200){\vector(0,1){15}} \put(65,40){\vector(-1,0){15}}
\put(10,25){\line(1,1){295}} \put(5,320){$k_2$} \put(305,20){$k_1$}
\put(70,20){$\beta$} \put(150,20){$-\alpha-1$} \put(0,125){$\beta+1$}
\put(0,205){$-\alpha$}

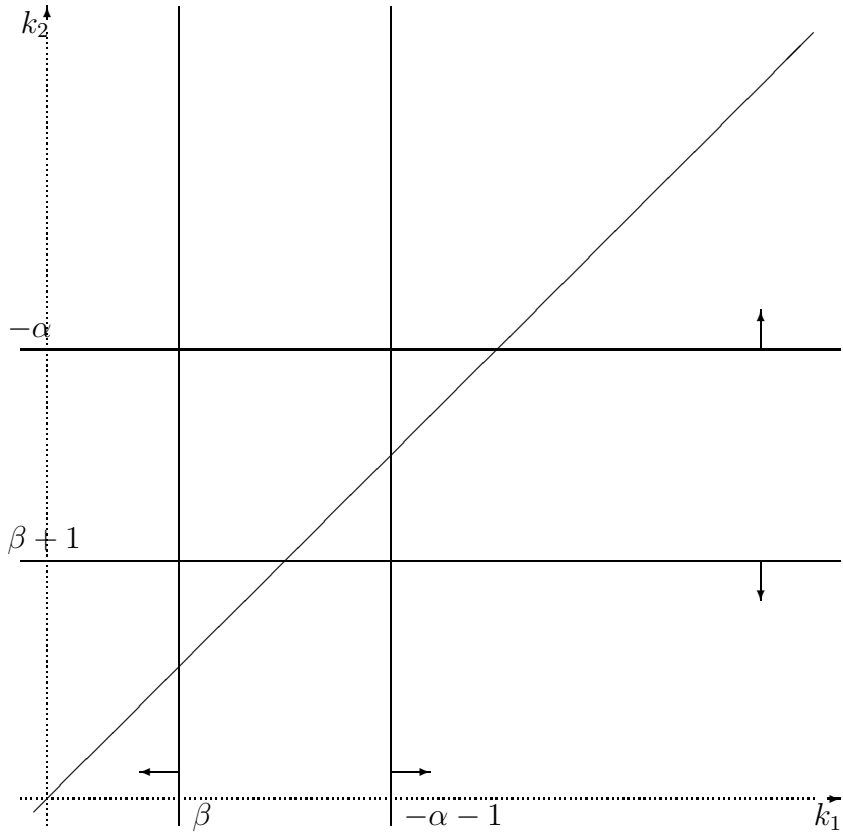
\put(100,0){Fig.3. The structure of $\pi_{\alpha,\beta}$ with
$\alpha+\beta=-2$.}
\end{picture}

\medskip

Case 4. $\alpha+\beta \leq -3$. In this case the intersection point of
$\mathcal{L}_1^-$ and $\mathcal{L}_2^+$ belongs to $\mathbf K^+$ (see
Fig.4). Also, there are simple submodules $V^s_1$, $V^s_2$, $V^s_3$ in
$\mathbb C[S(\mathbb D)]_q$, but $\mathbb C[S(\mathbb D)]_q$ does not
decompose into their direct sum.

\begin{picture} (300,330)(0,0)
\put(15,20){\dottedline{2}(0,0)(0,300)} \put(15,325){\vector(0,1){5}}
\put(5,30){\dottedline{2}(0,0)(300,0)} \put(305,30){\vector(1,0){5}}
\put(5,100){\line(1,0){310}} \put(5,180){\line(1,0){310}}
\put(65,20){\line(0,1){310}} \put(125,20){\line(0,1){310}}
\put(125,40){\vector(1,0){15}} \put(300,100){\vector(0,-1){15}}
\put(300,180){\vector(0,1){15}} \put(65,40){\vector(-1,0){15}}
\put(10,25){\line(1,1){295}} \put(5,320){$k_2$} \put(305,20){$k_1$}
\put(70,20){$\beta$} \put(130,20){$-\alpha-1$} \put(0,105){$\beta+1$}
\put(0,185){$-\alpha$}

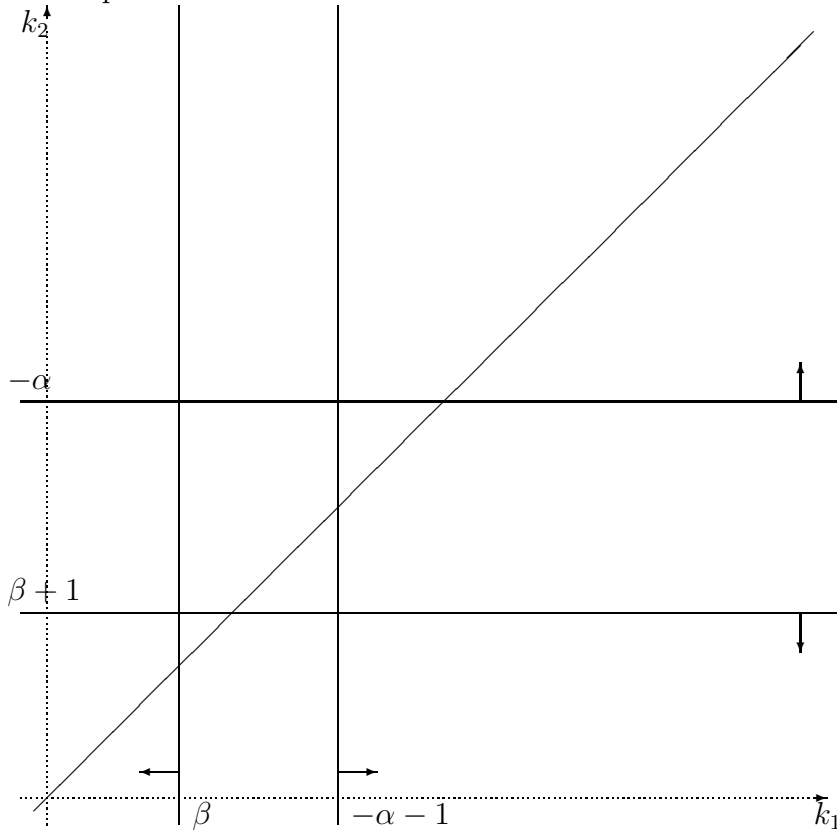
\put(100,0){Fig.4. The structure of $\pi_{\alpha,\beta}$ with
$\alpha+\beta \leq -3$.}
\end{picture}

\bigskip

Turn now to the general case. Consider all possible values of
$\alpha+\beta +n-1$.

{\bf Case 1.} $\alpha+\beta +n-1\geq 1$. In this case the hyperplanes
$\mathcal{L}_j^{\pm},\,j=1,\ldots,n$, bound in $\mathbf K^+$ a subset that
corresponds to a unique simple {\it finite dimensional} submodule
$$V^s=\bigoplus \limits_{\{\overline{\mathbf k} \in \widehat{K}|-\alpha-n+j
\leq k_j \leq \beta+j-1 \text{ for all } j=1, \ldots , n\}}
V_{\overline{\mathbf k}}.$$

{\bf Case 2.} $\alpha+\beta+n-1=0$. In this case the hyperplanes
$\mathcal{L}_j^{+}$ and $\mathcal{L}_j^{-}$ coincide. There are $n$ simple
submodules in $\mathbb C[S(\mathbb D)]_q$:
\begin{equation}\label{c_2_submod}
V^s_j=\bigoplus \limits_{\{\overline{\mathbf k} \in
\widehat{K}|k_j=\beta+j-1\}} V_{\overline{\mathbf k}}, \quad j=1,\ldots,n.
\end{equation}

{\bf Case 3.} $\alpha+\beta=-n$. Here the distance between
$\mathcal{L}_j^{+}$ and $\mathcal{L}_j^{-}$ is 1. This allows one to
decompose the set $\widehat{K}$ into a direct sum of $n+1$ subsets
$\widehat{K}_i,i=1,\ldots,n+1$, those correspond to the simple submodules:
$V^s_i= \bigoplus \limits_{\{\overline{\mathbf k} \in \widehat{K}_i\}}
V_{\overline{\mathbf k}} \;\subset \mathbb C[S(\mathbb D)]_q$. The subsets
$\widehat{K}_i$ are defined as follows:
$$
\widehat{K}_i=\{\overline{\mathbf k} \in \widehat{K}|k_{i-1} \geq
-\alpha-n +i-1, \beta+i-1 \geq k_i\}
$$
(for $i=1$ and $i=n+1$ we put respectively
$\widehat{K}_1=\{\overline{\mathbf k} \in \widehat{K}|k_1 \leq \beta\}$
and $\widehat{K}_{n+1}=\{\overline{\mathbf k} \in \widehat{K}|k_n \geq
-\alpha\}$).

\medskip

{\bf Case 4.} $\alpha+\beta+n-1 \leq -2$. Also, there are simple
submodules corresponded to $\widehat{K}_i$. However, $\mathbb C[S(\mathbb
D)]_q$ is not their direct sum.

Thus we have proved the following
\begin{corollary}\label{cor_4}
For $\alpha,\beta \in \mathbb Z$ only one of the representations
$\pi_{\alpha,\beta}$ and $\pi_{-n-\beta,-n-\alpha}$ has an irreducible
finite dimensional subrepresentation.
\end{corollary}

\subsection {Intertwining operators} \label{eqv}

In this section we construct the intertwining operators between the
representations $\pi_{\alpha,\beta}$ and $\pi_{-n-\beta,-n-\alpha}$ for
non-integral $\alpha,\beta$. This allows one to prove Proposition
\ref{pr_eqv}.

Let $A:\mathbb C[S(\mathbb D)]_q \rightarrow \mathbb C[S(\mathbb D)]_q$ be
an intertwining operator, i.e., for all $\xi \in U_q \mathfrak{sl}_{2n},\,
v \in \mathbb C[S(\mathbb D)]_q$, we have $A \pi_{\alpha,\beta}(\xi)(v)
=\pi_{-n-\beta,-n-\alpha}(\xi)(Av).$ The operators $\pi_{\alpha,\beta}(U_q
\mathfrak{k}_{ss})$ are independent of $\alpha,\beta$ and
$\pi_{\alpha,\beta}(K_n)=\pi_{-n-\beta,-n-\alpha}(K_n)$. Also,
$V_{\overline{\mathbf k}}$ and $V_{\overline{\mathbf m}}$ are
non-isomorphic $U_q \mathfrak{k}$-modules for $\overline{\mathbf k} \ne
\overline{\mathbf m}$. Then $A(\alpha,\beta)|_{V_{\overline{\mathbf k}}} =
a_{\overline{\mathbf k}} (\alpha,\beta)$, $a_{\overline{\mathbf k}}
(\alpha,\beta) \in \mathbb{C}$.  By the additional assumption
$a_{\overline {\mathbf 0}} (\alpha,\beta)=1$, we have the explicit
formulas for the coefficients $a_{\overline{\mathbf
k}}(\alpha,\beta)=A(\alpha,\beta)|_{V_{\overline{\mathbf k}}}$ of the
intertwining operator $A$
\begin{equation}\label{inter_op}
a_{\overline{\mathbf k}}(\alpha,\beta) = \prod_{j=1}^n P_j(\alpha,\beta),
\end{equation}
where \begin{equation*} P_j(\alpha,\beta) =
\begin{cases}
\prod
\limits_{i=0}^{k_j-1}\frac{1-q^{2(\alpha+n+i-j+1)}}{1-q^{2(-\beta+i-j+1)}},
&\text{for} \quad k_j >0,
\\ 1, & \text{for} \quad k_j=0,
\\ \prod \limits_{i=1+k_j}^{0} \frac{1-q^{2(-\beta+i-j)}}{1-q^{2(\alpha+n+i-j)}}, &
\text{for} \quad k_j < 0.
\end{cases}
\end{equation*}

For fixed $\alpha-\beta \in \mathbb Z$ $A$ is a meromorphic
operator-valued function with simple poles at integral points. Note that
$A$ coincides up to a multiplicative constant with the so called standard
intertwining operator, and the constant can be expressed from a $q$-analog
of the Harish-Chandra $c$-function.

\subsection{Unitarizable representations of the degenerate principal
series}\label{unit}

In this section we list necessary and sufficient conditions for modules of
the degenerate principal series and their simple submodules to be
unitarizable.

Recall the definition of unitarizable module. Let A be a $*$-Hopf algebra,
$W$ an $A$-module. An $A$-module $W$ is {\it unitarizable} if there exists
a scalar product\footnote{I.e., sesquilinear Hermitian-symmetric positive
form.} $(\cdot,\cdot)$, which is $A$-invariant, i.e.,
\begin{equation*}
(a u,v)=(u,a^{*}v) \quad \text{for any} \quad u,v \in W,\, a \in A.
\end{equation*}

\medskip

We can present the following series of simple unitary {\it representations
of the degenerate principal series} related to the Shilov boundary.

\medskip

{\it The principal unitary series:} $\mathrm{Re} (\alpha+\beta)=-n$,
$\alpha,\beta \not \in \mathbb Z$. In this case all representations are
unitarizable. The invariant scalar product is obtained from the results of
\cite{Cauhy_Sege}.

{\it The complementary series:} $\mathrm{Im} (\alpha+\beta)=0$,
$[-\alpha-n]=[\beta]$, $\alpha,\beta \not \in \mathbb Z$. In this case the
representations $\pi_{\alpha,\beta}$ are unitarizable too.

{\it The strange series:} $\mathrm{Im} \alpha=\frac{\pi}{h}$. For such
values of the parameters the respective representations
$\pi_{\alpha,\beta}$ are irreducible and unitarizable. This series of
representations has no classical analog (cf. \cite{KL_Pak}).

\medskip

Now let $\alpha,\beta \in \mathbb{Z}$. (Recall that in this case
$\pi_{\alpha,\beta}$ is reducible.) For such $\alpha,\beta$ there might
exist unitarizable simple submodules in the respective module (we will
mention them below), although the module is not unitarizable. Consider all
possible cases:

{\bf Case 1.} $\alpha+\beta\geq 2-n$. In this case the representation is
not unitarizable and its unique irreducible subrepresentation is not
unitarizable too.

{\bf Case 2.} $\alpha+\beta=1-n$. In this case there exist $n$ irreducible
unitarizable subrepresentations of the representation
$\pi_{\alpha,1-n-\alpha}$. Precisely, $V^s_j$ (see (\ref{c_2_submod})) is
a simple submodule in $\mathbb C[S(\mathbb D)]_q$ for any $j=1,\ldots,n$.
Notice that each $V^s_j$ can be equipped with a $U_q
\mathfrak{su}_{n,n}$-invariant scalar product $(\cdot,\cdot)$. Such
modules are called small representations because they have ''poor''
decompositions into isotypic components.

{\bf Case 3.} $\alpha+\beta=-n$. In this case the representations are
completely reducible, their irreducible subrepresentations
$V^s_i,\,i=1,\ldots,n+1$, (see Section \ref{reduce}) are unitarizable
(actually, the required invariant scalar product is the same as that for
the principal unitary series).

{\bf Case 4.} $\alpha+\beta \leq -1-n$. In this case the submodules
$V^s_i,\,i=1,\ldots,n+1$, are unitarizable, although there exist
non-unitarizable quotient modules in $\mathbb C[S(\mathbb D)]_q$.

\section*{Acknowledgment}
The author would like to thank L.L. Vaksman for constant attention and
helpful discussions.

\end{document}